\documentclass[11pt,times,letter]{article}
\usepackage{amsmath, amsfonts,amsthm,amssymb,bm}
\usepackage{dsfont}

\usepackage[english]{babel}

\usepackage{latexsym,amssymb,amsfonts,graphicx}
\usepackage{amsmath,dsfont}
\usepackage{verbatim}
\usepackage{mathrsfs}
\usepackage{bm}
\usepackage{color}
\usepackage{epsfig}

\usepackage{url}

\usepackage{bm}
\usepackage[colorlinks,linkcolor=blue,citecolor=blue,anchorcolor=blue]{hyperref}

\newcommand{\E}{\mathbb{E}}
\newcommand{\PP}{\mathbb{P}}
\newcommand{\RR}{\mathbb{R}}

\newtheorem{theorem}{Theorem}[section]
\newtheorem{lemma}{Lemma}[section]
\newtheorem{coro}{Corollary}[section]

\newtheorem{prop}{Proposition}[section]
\newtheorem{rem}{Remark}[section]

\def\F{{\cal F}}
   \def\D{d}     
 \def\N{\nabla}

 \def\2{\oplus} \def\3{\otimes} \def\4{\ominus}
\def\5{\circ} \def\6{\odot} \def\7{\backslash} \def\8{\infty}
\def\9{\bigcap} \def\0{\bigcup} \def\+{\pm} \def\-{\mp}
  \def\0{\bf 0}

\newcommand{\1}{{\mbox{$1$}\rm \!l}}
\allowdisplaybreaks[4]

\definecolor{Red}{rgb}{1.00, 0.00, 0.00}

\definecolor{DRed}{rgb}{0.5, 0.00, 0.00}

\definecolor{Blue}{rgb}{0.00, 0.00, 1.00}

\definecolor{Green}{rgb}{0.0, 0.4, 0.0}

\title
{
Stochastic Delay Differential Equations with Jump Reflection: Invariant Measure
}

\author{
Lijun Bo \thanks{Email: lijunbo@ustc.edu.cn, School of Mathematical Sciences, University of Science and Technology of China, Hefei, Anhui
Province, 230026, China.} \and
Chenggui Yuan\thanks{Email: C.Yuan@swansea.ac.uk, Department of Mathematics, Swansea University,
Swansea SA2 8PP, UK.}}

\begin{document}
\maketitle
\begin{abstract}
In this paper, we consider a class of multi-dimensional stochastic delay differential equations with jump reflection. Based on existence and uniqueness of the strong solution to the equation, we prove that the Markov semigroup generated by the segment process corresponding to the solution admits a unique invariant measure on the Skorohod space when the coefficients of equation satisfy a class of monotone conditions. Finally, we establish a relationship between the regulator and the local time of the solution and discuss a local time property at large time under the stationary setting.\\


\end{abstract}

\section{Introduction}\label{sec:intro}
We take  a complete filtered probability space $(\Omega,\F,{\mathbb F}=(\F_t;\ t\geq0),\PP)$ carrying an ${\mathbb F}$-adapt
$n\in{\mathbb{N}}$-dimensional Brownian motion $W=(W^l(t);\ t\geq0)_{l=1,\ldots,n}$ and an ${\mathbb F}$-adapt Poisson random
measure $(N((0,t]\times A);\ A\in{\mathcal{B}}({\mathcal{E}}))$ with
intensity $(\nu(A)t;\ A\in{\mathcal{B}}({\mathcal{E}}))$, where
$t>0$ and ${\mathcal{E}}=\RR^d/\{{\bf0}\}$. Here the filtration
${\mathbb F}$ is assumed to satisfy the usual conditions. For $I\subset \RR$, let $D(I; \RR_+^d)$ denote the space of all right-continuous
functions with left limits (r.c.l.l.) from $I$ to $\RR_+^d$ with Skorohod topology, see also Liptser and Shiryayev~\cite{LS}. In this paper, we consider the following $d\in{\mathbb{N}}$-dimensional reflected stochastic delay differential equation with jumps:
\begin{align}\label{eq:rsddej}
\begin{cases}
\D X(t)= b(t,X(t),X(t-\tau))\D t + \sigma(t,X(t),X(t-\tau))\D
W(t)\\
\qquad\qquad+\int_{\mathcal{E}}g(t,X(t-),X((t-\tau)-),\rho)\tilde{N}(\D\rho,\D
t)+\D
K(t), & {\rm on}\ t\geq0,\\
X(t) = \xi(t)\in\RR_+^d, & {\rm on}\ t\in[-\tau,0],
\end{cases}
\end{align}
where $\tau>0$ is a deterministic delay level, the initial data $\xi(\cdot)\in D([-\tau,0];\RR_+^d)$, the
coefficients $(b(\cdot),g(\cdot))$ are the functions of $d$-dimensioanl column vectors and the volatility $\sigma(\cdot)$ is the function of
the $d\times n$-matrix. $\tilde{N}(\D\rho,\D t):={N}(\D\rho,\D t)-\nu(\D\rho)\D t$ defines the compensated version of the Poisson measure ${N}(\D\rho,\D t)$ where characteristic measure $\nu(\D\rho)$ is a $\sigma$-finite measure on $({\mathcal{E}},{\mathcal{B}}({\mathcal{E}}))$. Here $K=(K^i(t);\ t\geq0)_{i=1,\ldots,d}$ is a $d$-dimensional
nonnegative process, which is called the regulator for the
$d$-dimensional solution process $X=(X(t);\ t\geq-\tau)$ at the orthant. Moreover, the
regulator $K$ can be uniquely determined by the following properties
up to a positive constant factor (see Kinnally and Williams~\cite{KW}):
\begin{itemize}
  \item[{\sf(a)}] For $i=1,\dots,d$, the paths of $t\to
K^i(t)$ are non-decreasing, r.c.l.l. and $K^i(0)=K^i(0-)={0}$;
\item[{\sf(b)}] For all $t\geq0$, it holds that
\begin{eqnarray}\label{eq:Ksupport}
\int_0^t\left<X(s),\D K(s)\right>=0,
\end{eqnarray}
where $\left<x,y\right>=\sum_{i=1}^dx_iy_i$ for $x=(x_i)_{d\times1}$
and $y=(y_i)_{d\times1}\in\RR^d$.
\end{itemize}

For any $d\times n$-matrix $a=(a_{ij})_{d\times n}$, define
$\|a\|=\sqrt{{\rm Tr}[aa^{\top}]}$, where $a^{\top}$ is the
transpose of $a$ and ${\rm Tr}[aa^{\top}]$ denotes the trace of the
matrix $aa^{\top}$. Define $|x|=\sqrt{\left<x,x\right>}$ for any
$x=(x_i)_{d\times1}\in\RR^d$. We work in the following assumptions
on coefficient functions $(b(\cdot),\sigma(\cdot),g(\cdot))$ in Eq.~\eqref{eq:rsddej} throughout the paper:
\begin{enumerate}
 \item[{\sf (A1)}] there exists a constant $\alpha>0$ and constants $\alpha_1>\alpha_2>0$ such that
 \begin{eqnarray*}
2\langle x, b(t, x, y) \rangle
+\|\sigma(t,x,y)\|^2+\int_{\mathcal{E}}|g(t,x,y,\rho)|^2\nu(\D\rho)\le
\alpha-\alpha_1|x|^2+\alpha_2|y|^2,
\end{eqnarray*}
for all $t\in\RR_+$ and $x,y\in\RR_+^d$.
 \item[{\sf (A2)}]  there exist two
constants $\beta_1>\beta_2>0$ such that
\begin{eqnarray*}\label{eq:assumption}
&&2\left<x-\hat{x},b(t,x,y)-b(t,\hat{x},\hat{y})\right>+\|\sigma(t,x,y)-\sigma(t,\hat{x},\hat{y})\|^2\nonumber\\
&&\quad
+\int_{\mathcal{E}}|g(t,x,y,\rho)-g(t,\hat{x},\hat{y},\rho)|^2\nu(\D\rho)\leq-\beta_1|x-\hat{x}|^2+\beta_2|y-\hat{y}|^2,
\end{eqnarray*}
for all $t\in\RR_+$ and $x,\hat{x},y,\hat{y}\in\RR_+^d$.
\end{enumerate}
An illustrative example for conditions {\sf(A1)} and {\sf(A2)} is
to take the drift coefficient $b(t,x,y)=-\gamma(t) x+ \theta(t) y$
with $\gamma(t),\theta(t)>0$. We assume that
$\gamma_*=\inf_{t\geq0}\gamma(t)$ and
$\theta^*=\sup_{t\geq0}\theta(t)$ are finite. For
$(t,x,\hat{x},y,\hat{y})\in\RR_+\times\RR_+^d\times\RR_+^d\times\RR_+^d\times\RR_+^d$,
it holds that
\begin{eqnarray*}
2\left<x-\hat{x},b(t,x,y)-b(t,\hat{x},\hat{y})\right>\leq
-(2\gamma_*-\varepsilon^2)|x-\hat{x}|^2+\frac{|\theta^*|^2}{\varepsilon^2}|y-\hat{y}|^2,
\end{eqnarray*}
for any $\varepsilon>0$. Further the coefficient
$(\sigma(\cdot),g(\cdot))$ are assumed to satisfy $\sigma(t,0,0)=0$, $g(t,0,0,\rho)=0$ for all $(t,\rho)\in\RR_+\times{\mathcal{E}}$
and the following Lipschitzian-type conditions:
\begin{eqnarray*}
\|\sigma(t,x,y)-\sigma(\hat{x},\hat{y},\rho)\|^2&\leq&\ell_\sigma(t)(|x-\hat{x}|^2+|y-\hat{y}|^2),\nonumber\\
|g(t,x,y,\rho)-g(t,\hat{x},\hat{y},\rho)|^2&\leq&
\ell_{g}(t,\rho)(|x-\hat{x}|^2+|y-\hat{y}|^2),
\end{eqnarray*}
where $\ell_\sigma(t)$ and $\ell_g(t,\rho)$ are positive functions
satisfying
$\ell_{\sigma,g}^*:=\sup_{t\geq0}(\ell_{\sigma}(t)+\int_{{\mathcal{E}}}\ell_g(t,\rho)\nu(\D\rho))<+\infty$.
If $\gamma_*-\ell_{\sigma,g}^*>\theta^*$, we can always find a
constant $\varepsilon>0$ such that
$\alpha_1:=2\gamma_*-\varepsilon^2-\ell_{\sigma,g}^*>\alpha_2:=\frac{|\theta^*|^2}{\varepsilon^2}+\ell_{\sigma,g}^*$.
Thus assumptions {\sf(A1)} and {\sf(A2)} are satisfied.

For the given Brownian motion $W$ and Poisson random measure $N$, we
call the ${\mathbb F}$-pair of  r.c.l.l. processes
$(X,K)=((X(t);\ t\geq-\tau),(K(t);\ t\geq0))$ is a strong solution
to Eq.~\eqref{eq:rsddej}, if they solve the following stochastic
integral equation:
\begin{equation}\label{ity}
\begin{cases}
 X(t) = \xi(0)+\int_0^tb(s,X(s),X(s-\tau))\D s + \int_0^t \sigma(s,X(s),X(s-\tau))\D
W(s)\\
\qquad\qquad+\int_0^t\int_{\mathcal{E}}g(s,X(s-),X((s-\tau)-),\rho)\tilde{N}(\D\rho,\D
s)+K(t)\in\RR_+^d, & {\rm on}\ t\geq0,\\
X(t)= \xi(t)\in\RR_+^d, & {\rm on}\ t\in[-\tau,0],
\end{cases}
\end{equation}
and the $\RR_+^d$-valued regulator $K$ satisfies properties {\sf
(a)} and {\sf(b)}. We first have the following remark on existence and
uniqueness of the strong solution to Eq.~\eqref{eq:rsddej}.
\begin{rem}
Under assumptions {\sf(A1)} and {\sf(A2)}, Eq.~\eqref{eq:rsddej} admits a unique strong solution defined as above.
As a matter of fact, the existence of the unique strong solution to
Eq.~\eqref{eq:rsddej} can be guaranteed by the following weaker
conditions than {\sf(A1)} and {\sf(A2)}, namely
\begin{enumerate}
 \item[{\sf (A1')}] there exists a constant $\bar \alpha>0$ such that
 \begin{eqnarray*}
2\langle x, b(t, x, y) \rangle
+\|\sigma(t,x,y)\|^2+\int_{\mathcal{E}}|g(t,x,y,\rho)|^2\nu(\D\rho)\le
\bar \alpha(1+|x|^2+|y|^2),
\end{eqnarray*}
for all $t\in\RR_+$ and $x,y\in\RR_+^d$.
 \item[{\sf (A2')}]  there exist two
constants $\bar \beta_1, \bar \beta_2>0$ such that
\begin{eqnarray*}\label{eq:assumption}
&&2\left<x-\hat{x},b(t,x,y)-b(t,\hat{x},\hat{y})\right>+\|\sigma(t,x,y)-\sigma(t,\hat{x},\hat{y})\|^2\nonumber\\
&&\quad
+\int_{\mathcal{E}}|g(t,x,y,\rho)-g(t,\hat{x},\hat{y},\rho)|^2\nu(\D\rho)\leq
\bar \beta_1|x-\hat{x}|^2+\bar \beta_2|y-\hat{y}|^2,
\end{eqnarray*}
for all $t\in\RR_+$ and $x,\hat{x},y,\hat{y}\in\RR_+^d$.
\end{enumerate}
The proof is similar to that of Krylov~\cite{K99} and von Renesse and Scheutzow~\cite{MM}, we omit it
here. Our aim is to use assumptions {\sf(A1)} and {\sf(A2)} to study invariant measures for Eq.~\eqref{eq:rsddej}.
\end{rem}

For $i=1,\dots,d$, let $K^{i,c}(t)=K^{i}(t)-\sum_{s\leq t}\Delta
K^i(s)$ be the continuous part of the $i$th-regulator $K^i(t)$, where
the $t$-time jump's size $\Delta K^i(t)=K^i(t)-K^i(t-)$ with left
limit $K^i(t-):=\lim_{s\uparrow t}K^i(s)$. It will be seen that the
continuous counterpart $K^{i,c}=(K^{i,c}(t);\ t\geq0)$ behaves like the
local time of the $i$th-element $X^i$ of the solution process $X$
when $X^i$ is treated as a r.c.l.l. semimartingale (see
Section~\ref{sec:local-time} below). However, the
jump of $K^i$ happens when $X^i-K^i$ jumps down below barrier zero
due to the appearance of some negative jump. This phenomenon is
usually called ``jump reflection" in the literature (see e.g.,
Slomi\'nski and Wojciechowski \cite{sw10} and Nam \cite{Nam10}).
Moreover, the corresponding jump's size of the $i$th-regulator is
given by
\begin{eqnarray}\label{eq:jumpK}
\Delta K^i(t) =
\left[\int_{\mathcal{E}}g_i(t,X(t-),X((t-\tau)-),\rho)N(\D\rho,\{t\}) + X^i(t-)\right]^-,
\end{eqnarray}
where $[x_i]^-=\max\{-x_i,0\}$ for $x_i\in\RR$ and
$[x]^-=(x_i^-)_{d\times1}$ for $x\in\RR^d$. Write $K^c=(K^{i,c})_{i=1,\ldots,d}$. From the ``jump reflection",
for all $t\geq0$, it also holds that
\begin{eqnarray}\label{eq:Kcsupport}
\int_0^t\left<X(s),\D K^c(s)\right>=0.
\end{eqnarray}

The similar one-sided Lipschitzian condition and monotone condition
as assumptions {\sf(A1)} and {\sf(A2)} have been discussed in Bao et al.
\cite{BTY} for stochastic delay equation without jump reflection, in
Marin-Rubio and Real \cite{MRR} and Zhang \cite{Z} for reflected stochastic differential equation without jumps. In particular, the
Picard's successive approximation used in Xu and Zhang \cite{XuZh}
can deal with existence and uniqueness of strong solutions to
Eq.~\eqref{eq:rsddej} when the regulator $K$ is described as a
local time. There exist negative jumps in Eq.~\eqref{eq:rsddej},
the regulator $K$ has jumps whose sizes can be
identified by Eq.~\eqref{eq:jumpK}. Let $(X_n;\
n=\{0\}\cup{\mathbb{N}})$ be the corresponding Picard's
approximating sequence to the $d$-dimensional solution processes.
Then the successive approximation to the jump $\Delta K(t)$ can be
established via \eqref{eq:jumpK}, namely
\[
\Delta
K_{n}(t)=\left[\int_{\mathcal{E}}g(t,X_{n-1}(t-),X_{n-1}((t-\tau)-),\rho)N(\D\rho,\{t\}) + X_{n-1}(t-)\right]^-,\ \ \ n\in{\mathbb{N}}.
\]
See Proposition 2.4 in Slomi\'nski and Wojciechowski \cite{sw10} for
more details. The literature on stochastic delay equations with (or without) jumps is extensive (see
e.g., \cite{AS10}, \cite{FO11}, \cite{KW}, \cite{OS10}, \cite{OS07},
\cite{HMS}, \cite{RRG}, \cite{WA11} and references therein).
Recently Kinnally and Williams \cite{KW} discussed existence and
uniqueness of stationary solutions to a class of reflected stochastic differential delay equation (SDDE)
driven by Brownian motions. However the stability in distribution implies the
existence and uniqueness of invariant measures for the corresponding
segment processes. To the best of our knowledge, it seems that there
exists not much literature to investigate SDDE with jump reflection.

An outline of the paper is as follows. Section~\ref{sec:moment} establishes an
estimate for the second-order moment associated to the segment
process of the solution to Eq.~\eqref{eq:rsddej} and then considers the
exponential integrability of the solution. Existence and uniqueness of invariant measures associated with $d$-dimensional
segment process is proved in Section~\ref{sec:invariant-measure}. In Section~\ref{sec:local-time}, we discuss the
relationship between the regulator $K$ and the local time of the solution process and the local time property in the stationary setting.

{\bf Additional Notation.}\quad For $f\in D(I;\RR_+^d)$ and $I\subset\RR$, $\|f\|_{I}:=\sup_{t\in I}|f(t)|$.
For the $d$-dimensional solution process $X=(X(t);\ t\geq-\tau)$,
the corresponding segment process $(X_t;\ t\geq0)$ is defined as $X_t(\theta)=X(t+\theta)$ with $\theta\in[-\tau,0]$, correspondingly
$\|X_t\|_{[-\tau,0]}:=\sup_{-\tau\leq\theta\leq0}|X(t+\theta)|$. Throughout the
paper, we use the conventions $\int_c^d:=\int_{(c,d]}$ and $\int_{c}^\infty:=\int_{(c,\infty)}$ for any real numbers $c<d$.

\section{Moment Estimates of Segment Process}
\label{sec:moment}

This section concentrates on the estimates of the second-order moment of the segment process $(X_t;\ t\geq0)$ and the
exponential moment for the solution process $(X(t);\ t\geq-\tau)$ to Eq.~\eqref{eq:rsddej}.

Before presenting these moment estimates, we first present the following auxiliary results which will serve to establish final estimates.

\begin{lemma}\label{lemma:ito-formula}
Let $X=(X(t);\ t\geq-\tau)$ be the strong solution to Eq.~\eqref{eq:rsddej}. Then, for any $F\in C^2(\RR_+^d)$ and $t\geq0$, it holds that
\begin{align}\label{eq:ito-formula2}
F(X(t))&=F(\xi(0)) + \int_0^t\left<\nabla
F(X(s)),b(s,X(s),X(s-\tau))\right>\D s\nonumber\\
&\quad+ \int_0^t\left<\nabla F(X(s)),\sigma(s,X(s),X(s-\tau))\D
W(s)\right>+\int_0^t\left<\nabla F(X(s)),\D K^c(s)\right>\nonumber\\
&\quad+ \int_0^t\int_{\mathcal{E}}\left<\nabla
F(X(s)),g(s,X(s),X(s-\tau),\rho)\right>\tilde{N}(\D\rho,\D s)\nonumber\\
&\quad+\frac{1}{2}\int_0^t {\rm Tr}\left[(\sigma\sigma^{\top})(s,X(s),X(s-\tau))D^2F(X(s))\right]\D s\nonumber\\
&\quad+\int_0^t\int_{\mathcal{E}}[F([X(s-)+g(s,X(s-),X((s-\tau)-),\rho)]^+)-F(X(s-))\nonumber\\
&\quad -\left<\nabla
F(X(s-)),g(s,X(s-),X((s-\tau)-),\rho)\right>]N(\D\rho,\D s),
\end{align}
where $[x]^+=(x_i^+)_{d\times1}$ and $x_i^+=\max\{x_i,0\}$ with
$x=(x_i)_{d\times1}\in\RR^d$, $\N F(x)$ denotes the gradient of
$F(x)$, $D^2F(x)$ is the the $d\times d$-matrix of second-order
partial derivatives of $F(x)$ and $K^c(t)$
corresponds to the continuous component of the regulator
$K(t)$ with $t\geq0$.
\end{lemma}

\noindent{\it Proof.}\quad By virtue of It\^{o} formula with jumps
(see e.g., Theorem IV.48 in Protter \cite{Protter}, Page 193), we have for $t\geq0$,
\begin{align}\label{eq:ito-formula}
F(X(t))&=F(\xi(0)) + \int_0^t\left<\nabla F(X(s-)),\D
X(s)\right>+\frac{1}{2}\int_0^t {\rm Tr}\left[(\sigma\sigma^{\top})(s,X(s),X(s-\tau))D^2F(X(s))\right]\D s\nonumber\\
&+\sum_{0<s\leq t}[F(X(s-)+\Delta X(s))-F(X(s-))-\left<\nabla
F(X(s-)),\Delta X(s)\right>].
\end{align}
For $i=1,\dots,d$, define the process with pure jumps:
\begin{eqnarray}\label{eq:Y-pure-jump}
Y^i(t)=\int_0^t\int_{\mathcal{E}}g_i(s,X(s-),X((s-\tau)-),\rho){N}(\D\rho,\D
s),\ \ \ \ \ t\geq0.
\end{eqnarray}
In terms of Eq.~\eqref{eq:rsddej}, the random jump amplitude of the
$i$th-element $X^i$ ($i=1,2,\dots,d$) is given by $\Delta
X^i(t)=\Delta Y^i(t) + \Delta K^i(t)$ for $t>0$. Using the following
representation of the jump's size of the $i$th-regulator $K^i$ (see also \eqref{eq:jumpK}), for $i=1,2,\dots,d$,
\begin{eqnarray}\label{eq:jump-Ki}
\Delta K^i(t) = [\Delta Y^i(t) + X^i(t-)]^-\ \ \ \ \ \ t>0,
\end{eqnarray}
we arrive at
\begin{eqnarray}\label{eq:jump-Xi}
\Delta X^i(t)=\Delta Y^i(t) + [\Delta Y^i(t) +
X^i(t-)]^-=: \varphi(X^i(t-),\Delta Y^i(t)),\ \ \ \ t>0.
\end{eqnarray}
The function $\varphi(x,y)=(\varphi(x_i,y_i))_{d\times1}$ with
$x=(x_i)_{d\times1}\in\RR_+^d$ and $y=(y_i)_{d\times1}\in\RR^d$,
where
$\varphi(x_i,y_i)=-x_i\1_{\{x_i+y_i\leq0\}}+y_i\1_{\{x_i+y_i>0\}}$
for $i=1,\dots,d$. Using the equality
$x_i+\varphi(x_i,y_i)=[x_i+y_i]^+$ and substitute the following
equality into \eqref{eq:ito-formula},
\begin{eqnarray*}
&&F(X(s-)+\Delta X(s))-F(X(s-))-\left<\nabla F(X(s-)),\Delta
X(s)\right>\nonumber\\
&=&F([X(s-)+\Delta Y(s)]^+)-F(X(s-))-\left<\nabla F(X(s-)),\Delta
Y(s)\right>-\left<\nabla F(X(s-)),\Delta K(s)\right>,
\end{eqnarray*}
we obtain \eqref{eq:ito-formula2}, where we have used
the finite variation property of the regulator $K$
and the following equality:
\begin{eqnarray*}
\int_0^t\left<\nabla F(X(s-)),\D K(s)\right>-\sum_{0<s\leq
t}\left<\nabla F(X(s-)),\Delta K(s)\right>=\int_0^t\left<\nabla
F(X(s)),\D K^c(s)\right>,\ \ \ \ t>0.
\end{eqnarray*}
Thus we complete the proof of the lemma. \hfill$\Box$

We further have the following corollary
\begin{coro}\label{coro:ito1}
Let $\lambda\in\RR$. Then the solution
process of Eq.~\eqref{eq:rsddej} admits for $t\geq0$,
\begin{align}\label{eq:ito1}
&e^{\lambda t}|X(t)|^2\leq|\xi(0)|^2 +
2\int_0^te^{\lambda s}\left<X(s),b(s,X(s),X(s-\tau))\right>\D s+ \lambda\int_0^te^{\lambda s}|X(s)|^2\D s\\
&\quad+\int_0^te^{\lambda s}\D M(s) + \int_0^te^{\lambda s} \|\sigma(s,X(s),X(s-\tau))\|^2\D
s+\int_0^t\int_{\mathcal{E}}e^{\lambda s}|g(s,X(s),X(s-\tau),\rho)|^2\nu(\D \rho)\D s,\nonumber
\end{align}
where the process $M=(M(t);\ t\geq0)$ is defined by
\begin{align*}
M(t)&:=\int_0^t\int_{\mathcal{E}}\big[\left|[X(s-)+g(s,X(s-),X((s-\tau)-),\rho)]^+\right|^2-|X(s-)|^2\big]\tilde{N}(\D\rho,\D
s)\nonumber\\
&\quad
+2\int_0^t\left<X(s),\sigma(s,X(s),X(s-\tau))\D
W(s)\right>.
\end{align*}
\end{coro}

\noindent{\it Proof.}\quad For $x\in\RR_+^d$, take $F(x)=|x|^2$ in Lemma \ref{lemma:ito-formula}. Then the equality \eqref{eq:ito-formula2} reads
\begin{align*}\label{eq:ito-formula3}
|X(t)|^2&=|\xi(0)|^2 + 2\int_0^t\left<X(s),b(s,X(s),X(s-\tau))\right>\D s
+ \int_0^t \|\sigma(s,X(s),X(s-\tau))\|^2\D s \nonumber\\
&\quad+ M(t)+\int_0^t\int_{\mathcal{E}}[\left|[X(s-)+g(s,X(s-),X((s-\tau)-),\rho)]^+\right|^2-|X(s-)|^2\nonumber\\
&\qquad\qquad\qquad\quad
-2\left<X(s-),g(s,X(s-),X((s-\tau)-),\rho)\right>]\nu(\D\rho)\D s,
\end{align*}
where we used the support property \eqref{eq:Kcsupport}, namely $\int_0^t\left<X(s),\D K^c(s)\right>=0$ for $t\geq0$. For $(t,x,y,\rho)\in\RR_+\times\RR_+^d\times\RR_+^{d}\times{\mathcal{E}}$, applying the following inequality
\begin{eqnarray}\label{ineq:1}
\left|[x+g(t,x,y,\rho)]^+\right|^2-|x|^2-2\left<x,g(t,x,y,\rho)\right>\leq
|g(t,x,y,\rho)|^2,
\end{eqnarray}
and we can conclude the validity of the inequality \eqref{eq:ito1} by using the
integration by parts. \hfill$\Box$

Corollary \ref{coro:ito1} can be used to establish the uniform estimate
of the second-order moment for the segment process
$(X_t;\ t\geq0)$.
\begin{prop}\label{prop:moment-segment}
Under the assumption {\sf (A1)}, it holds that
\begin{eqnarray}\label{eq:moment-segment}
\sup_{t\in[0,\infty)}\E\left[\|X_t\|_{[-\tau,0]}^2\right]<+\infty.
\end{eqnarray}
\end{prop}

\noindent{\it Proof.}\quad Let $\lambda>0$. By \eqref{eq:ito1} in Corollary \ref{coro:ito1}, we have under {\sf (A1)} that for $t\geq0$,
\begin{align*}
&\E\left[e^{\lambda t}|X(t)|^2\right]\leq\E\left[|\xi(0)|^2\right] +
2\E\left[\int_0^te^{\lambda s}\left<X(s),b(s,X(s),X(s-\tau))\right>\D s\right]\nonumber\\
&\qquad+ \E\left[\lambda\int_0^te^{\lambda s}|X(s)|^2\D s\right]
+\E\left[\int_0^te^{\lambda s} \|\sigma(s,X(s),X(s-\tau))\|^2\D s\right]\nonumber\\
&\qquad+\E\left[\int_0^t\int_{\mathcal{E}}e^{\lambda s}|g(s,X(s),X(s-\tau),\rho)|^2\nu(\D \rho)\D s\right]\nonumber\\
&\quad\leq\E\left[|\xi(0)|^2\right] + \alpha\int_0^te^{\lambda s}\D s+
(\lambda-\alpha_1)\E\left[\int_0^te^{\lambda s}|X(s)|^2\D s\right]+\alpha_2\E\left[\int_0^te^{\lambda s}|X(s-\tau)|^2\D s\right]\nonumber\\
&\quad=\E\left[|\xi(0)|^2\right] + \frac{\alpha}{\lambda}(e^{\lambda
t}-1)+e^{\lambda\tau}\E\left[\int_{-\tau}^0e^{\lambda v}|\xi(v)|^2\D
v\right]+\left(\lambda-\alpha_1+\alpha_2e^{\lambda\tau}\right)\E\left[\int_0^te^{\lambda
s}|X(s)|^2\D s\right].
\end{align*}
We next choose a constant $\lambda^*>0$ such that
$\lambda^*-\alpha_1+\alpha_2e^{\lambda^*\tau}=0$ since
$\alpha_1>\alpha_2$. Then for all $t\geq0$,
\begin{eqnarray}\label{yt1}
\E\left[|X(t)|^2\right]&\leq&e^{-\lambda^* t}\left\{\frac{\alpha}{\lambda^*}(e^{\lambda^*t}-1)+\E\left[|\xi(0)|^2\right]
+e^{\lambda^*\tau}\E\left[\int_{-\tau}^0e^{\lambda^* v}|\xi(v)|^2\D v\right]\right\}\nonumber\\
&\leq&\frac{\alpha}{\lambda^*}+\E\left[|\xi(0)|^2\right]
+e^{\lambda^*\tau}\E\left[\int_{-\tau}^0e^{\lambda^* v}|\xi(v)|^2\D v\right].
\end{eqnarray}
Let $\theta\in[-\tau,0]$. For any $t>\tau$, it follows from Lemma \ref{lemma:ito-formula} that
\begin{align*}
|X(t+\theta)|^2&=|X(t-\tau)|^2+
2\int_{t-\tau}^{t+\theta}\left<X(s),b(s,X(s),X(s-\tau))\right>\D s\nonumber\\
&\quad+\int_{t-\tau}^{t+\theta}\|\sigma(s,X(s),X(s-\tau))\|^2\D s+2\int_{t-\tau}^{t+\theta}\left<X(s),\sigma(s,X(s),X(s-\tau))\D W(s)\right>\nonumber\\
&\quad+2\int_{t-\tau}^{t+\theta}\int_{\mathcal{E}}\left<X(s),g(s,X(s),X(s-\tau),\rho)\right>\tilde{N}(\D\rho,\D s)\nonumber\\
&\quad+\int_{t-\tau}^{t+\theta}\int_{\mathcal{E}}[\left|[X(s-)+g(s,X(s-),X((s-\tau)-),\rho)]^+\right|^2-|X(s-)|^2\nonumber\\
&\qquad-2\left<X(s-),g(s,X(s-),X((s-\tau)-),\rho)\right>]N(\D\rho,\D s).
\end{align*}
Using the Burkh\"older-Davis-Gundy inequality, we obtain
\begin{align*}
&\E\left[\sup_{-\tau\leq\theta\leq0}\left|\int_{t-\tau}^{t+\theta}\left<X(s),\sigma(s,X(s),X(s-\tau))\D W(s)\right>\right|\right]\nonumber\\
&\qquad\leq\frac{1}{2}\E\left[\sup_{-\tau\leq\theta\leq0}|X(t+\theta)|^2\right]+C\E\left[\int_{t-\tau}^t\|\sigma(s,X(s),X(s-\tau))\|^2\D s\right],
\end{align*}
and
\begin{align*}
&\E\left[\sup_{-\tau\leq\theta\leq0}\left|\int_{t-\tau}^{t+\theta}
\int_{\mathcal{E}}\left<X(s),g(s,X(s-),X((s-\tau)-),\rho)\right>\tilde{N}(\D\rho,\D s)\right|\right]\nonumber\\
&\qquad\leq\frac{1}{4}\E\left[\sup_{-\tau\leq\theta\leq0}|X(t+\theta)|^2\right]+C\E\left[\int_{t-\tau}^t\int_{\mathcal{E}}|g(s,X(s),X(s-\tau),\rho)|^2\nu(\D\rho)\D
s\right],
\end{align*}
where $C>0$ is some positive constant. On the other hand, using \eqref{ineq:1}, it follows that
\begin{align*}
&\E\Bigg[\sup_{-\tau\leq\theta\leq0}\int_{t-\tau}^{t+\theta}\int_{\mathcal{E}}[|[X(s-)+g(s,X(s-),X((s-\tau)-),\rho)]^+|^2
-|X(s-)|^2\nonumber\\
&\quad\qquad-2\left<X(s-),g(s,X(s-),X((s-\tau)-),\rho)\right>]N(\D\rho,\D s)\Bigg]\nonumber\\
&\quad\leq\E\left[\sup_{-\tau\leq\theta\leq0}\int_{t-\tau}^{t+\theta}\int_{\mathcal{E}}|g(s,X(s-),X((s-\tau)-),\rho)|^2N(\D\rho,\D s)\right]\nonumber\\
&\quad\leq\E\left[\int_{t-\tau}^{t}\int_{\mathcal{E}}|g(s,X(s),X(s-\tau),\rho)|^2\nu(\D\rho)\D
s\right].
\end{align*}
Therefore it holds that
\begin{eqnarray}\label{eq:moment-single}
\E\left[\|X_t\|_{[-\tau,0]}^2\right]\leq4\E\left[|X(t-\tau)|^2\right]+C\int_{t-\tau}^t\E\left[|X(s-\tau)|^2\right]\D s,
\end{eqnarray}
for some positive constant $C$ which is independent of time $t$. The
required assertion follows from \eqref{yt1}.  \hfill$\Box$

The following result relates to the exponential moment of the solution process $X$, which is a reflected delay
version of R\"ochner and Zhang \cite{RZ07}'s exponential
integrability of the solution without reflection and delay when the
drift and diffusion coefficients $(b(\cdot),\sigma(\cdot))$ are uniformly bounded
on $(t,x,y)\in[0,T]\times\RR_+^d\times\RR_+^d$ with $T>0$. In
particular, we establish the exponential moment estimate of the
following reflected delay equation without drift and diffusive parts:
\begin{eqnarray}\label{eq:rsddej1}
X(t)&=&\xi(0) + \int_0^t\int_{\mathcal{E}} g(s,X(s-),X((s-\tau)-),\rho)\tilde{N}(\D\rho,\D s) + K(t)\nonumber\\
&= :&\xi(0) + \tilde{Y}(t) + K(t),\ \ \ \ \ t\geq0,
\end{eqnarray}
where the $d$-dimensional process $\tilde{Y}=(\tilde{Y}^i(t);\ t\geq0)_{i=1,\ldots,d}$ is the compensated version of the pure jump process $Y$ defined as the stochastic integral \eqref{eq:Y-pure-jump}.
\begin{lemma}\label{lemma:exp-integral}
For the characteristic measure and the jump coefficient  $(\nu(\cdot),g(\cdot))$, suppose that there exists a constant $\ell_g>0$ such that
\begin{eqnarray}\label{cond:g}
|g(t,x,y,\rho)|\leq \ell_gh(t,\rho),\ \ \ \ \ \ \forall\
(t,x,y,\rho)\in[0,T]\times\RR_+^d\times\RR_+^d\times{\mathcal{E}},
\end{eqnarray}
where $T>0$ and $h(t,\rho)$ is a nonnegative measurable function
satisfying
\begin{eqnarray}\label{cond:h}
\sup_{t\in[0,T]}\int_{{\mathcal{E}}}h^2(t,\rho)e^{\varrho
h(t,\rho)}\nu(\D\rho)<+\infty,
\end{eqnarray}
for $\varrho>0$. If drift and diffusion coefficients
$(b(\cdot),\sigma(\cdot))$ are uniformly bounded on
$(t,x,y)\in[0,T]\times\RR_+^d\times\RR_+^d$, then for any finite $a>0$, it holds that
\begin{eqnarray}\label{eq:expmoment}
\E\left[\exp\left(a\sup_{t\in[0,T]}|X(t)|\right)\right]<+\infty.
\end{eqnarray}
\end{lemma}

To prove Lemma~\ref{lemma:exp-integral}, we need the following lemma:
\begin{lemma}\label{coro:expIto}
Let $H\in C^2(\RR_+^d)$. Define the following function by
\begin{eqnarray}\label{eq:fcnQ}
Q^H(t,x,y)&:=&\left<\nabla H(x),b(t,x,y)\right>+\frac{1}{2}{\rm
Tr}\left[(\sigma\sigma^{\top})(t,x,y)(D^2H(x)+\nabla
H(x)\otimes\nabla H(x))\right]\\
&&+\int_{\mathcal{E}}\left[\exp\left\{{H([x+g(t,x,y,\rho)]^+)-H(x)}\right\}-1-\left<\nabla
H(x),g(t,x,y,\rho)\right>\right]\nu(\D\rho),\nonumber
\end{eqnarray}
on $(t,x,y)\in\RR_+\times\RR_+^d\times\RR_+^d$. We further define
the positive process by
\begin{eqnarray}\label{eq:exp-process}
Y^E(t):=\exp\left(H(X(t))-H(\xi(0))-\int_0^tQ^H(s,X(s),X(s-\tau))\D
s\right),\ \ \ \ t\geq0,
\end{eqnarray}
where $X=(X(t);\ t\geq-\tau)$ is the $d$-dimensional solution
process to Eq.~\eqref{eq:rsddej}. Then the positive process
$Y^E=(Y^E(t);\ t\geq0)$ satisfies
\begin{eqnarray}\label{eq:exp-martinagle}
Y^E(t)=1+M^E(t)+\int_0^t Y^E(s)\left<\nabla H(X(s)),\D
K^c(s)\right>,\ \ \ \ \ t\geq0.
\end{eqnarray}
Here the process $M^E=(M^E(t);\ t\geq0)$ is an ${\mathbb{F}}$-local martingale taking values on $\RR$.
\end{lemma}

\noindent{\it Proof.}\quad For $x\in\RR_+^d$, take the
function $F(x)=\exp(H(x))$ in Lemma \ref{lemma:ito-formula}. Since $\nabla F(x)=F(x)\nabla H(x)$ and
$D^2 F(x)=F(x)[D^2H(x)+\nabla H(x)\otimes\nabla H(x)]$,  we arrive
at
\begin{eqnarray*}
\frac{\D\exp\left\{H(X(t))-H(\xi(0))\right\}}{\exp\left\{H(X(t))-H(\xi(0))\right\}}
=Q^H(t,X(t),X(t-\tau))\D t + \D \hat{M}(t) + \left<\nabla H(X(t)),\D
K^c(t)\right>,
\end{eqnarray*}
where $\hat{M}=(\hat{M}(t);\ t\geq0)$ is a real-valued ${\mathbb{F}}$-local martingale. Then the equality
\eqref{eq:exp-martinagle} follows from applying the integration by parts
to
$\exp\left\{H(X(t))-H(\xi(0))\right\}\exp(-\int_0^tQ^H(s,X(s),X(s-\tau))\D
s)$. \hfill$\Box$

We are in the position to prove Lemma~\ref{lemma:exp-integral}.

\noindent{\it Proof of Lemma~\ref{lemma:exp-integral}.}\quad We adopt the test function used in R\"ochner and Zhang \cite{RZ07} to discuss our reflected delay
case. Consider the following function on $\RR_+^d$ given by
\begin{eqnarray}\label{eq:fcnH}
H_\lambda(x)=\sqrt{1+\lambda|x|^2},\ \ \ \ \ \lambda>0,\
x\in\RR_+^d.
\end{eqnarray}
Then the gradient $\nabla H_{\lambda}(x)=\lambda
H_\lambda^{-1}(x)x$ with $x\in\RR_+^d$. For any $i,j=1,\dots,d$,
and $x\in\RR_+^d$, the partial derivatives
\begin{eqnarray}\label{eq:derivaitve-esti}
\frac{\partial H_{\lambda}(x)}{\partial x_i}\leq\sqrt{\lambda},\ \ \ {\rm and}\ \
\frac{\partial^2 H_{\lambda}(x)}{\partial x_i\partial x_j}+\frac{\partial H_{\lambda}(x)}{\partial x_i}\frac{\partial H_{\lambda}(x)}{\partial x_j}
\leq2\lambda.
\end{eqnarray}
Recall the real-valued process $Y^E$ given by
\eqref{eq:exp-martinagle}. Note that, for all $t\geq0$,
\begin{eqnarray*}\label{eq:HKc}
0\leq\int_0^t\left<\nabla H_{\lambda}(X(s)),\D
K^c(s)\right>=\int_0^t\lambda H_\lambda^{-1}(X(s))\left<X(s),\D
K^c(s)\right>\leq\lambda\int_0^t\left<X(s),\D K^c(s)\right>=0,
\end{eqnarray*}
using the support property \eqref{eq:Kcsupport}. This implies
that
\begin{eqnarray*}
\int_0^t\left<\nabla H_{\lambda}(X(s)),\D K^c(s)\right>=0,\ \ \ \ \
\forall\ t\geq0.
\end{eqnarray*}
Then by Lemma \ref{coro:expIto}, we have that
\begin{eqnarray*}
Y_\lambda^E(t):=\exp\left(H_\lambda(X(t))-H_\lambda(\xi(0))-\int_0^tQ^{H_\lambda}(s,X(s),X(s-\tau))\D
s\right),\ \ \ \ t\geq0
\end{eqnarray*}
is a positive ${\mathbb F}$-local martingale and hence it is a
supermartingale, where the function $Q_\lambda^H(t,x,y)$ is given by
\eqref{eq:fcnQ} with the function $H(x)$ replaced by $H_\lambda(x)$.
Using the inequality $H_\lambda([x+g(t,x,y,\rho)]^+)\leq
H_\lambda(x+g(t,x,y,\rho))$ and the estimates
\eqref{eq:derivaitve-esti}, it is not difficult to prove that, for
any $(t,x,y)\in[0,T]\times\RR_+^d\times\RR_+^d$,
\begin{eqnarray*}
Q^{H_\lambda}(t,x,y)\leq C_1\left[1+\sup_{t\in[0,T]}\int_{\mathcal{E}}h^2(t,\rho)e^{\sqrt{\lambda}h(t,\rho)}\nu(\D\rho)\right]
:=C_2<+\infty,
\end{eqnarray*}
under conditions \eqref{cond:g} and \eqref{cond:h}, where
positive constants $C_1=C_1(d,\lambda)$ and $C_2=C_2(d,\lambda,T)$
depend on the dimension number $d$, the parameter $\lambda$ and the
time level $T$ only. Based on the above estimate of $Q^{H_\lambda}(t,x,y)$, the desired result follows from Proposition
4.2 and Corollary 4.3 in \cite{RZ07}.
\hfill$\Box$

\section{Invariant Measures}
\label{sec:invariant-measure}

This section will establish existence and uniqueness of an
invariant measure of the $d$-dimensional segment process $(X_t;\ t\geq0)$ under assumptions {\sf(A1)} and {\sf(A2)}.

For the initial data $\xi\in D([-\tau,0];\RR_+^d)$, we use
$(X_t^\xi;\ t\geq0)$ to represent the corresponding segment process
to the solution process $X=(X(t);\ t\geq-\tau)$ with
$X(t)=\xi(t)$ on $[-\tau,0]$. Define the Markov semigroup associated
with the segment process by
\begin{eqnarray}\label{eq:semigroup}
{\mathcal{P}}_tf(\xi)={\E}\big[f(X_t^\xi)\big],
\end{eqnarray}
for all bounded continuous functions $f$ defined on Skorohod space
$D([-\tau,0];\RR_+^d)$ with ${\mathcal{M}}$ being the Borel $\sigma$-algebra generated
by $D([-\tau,0];\RR_+^d)$ (see \cite{LS}). For any finite time level $T>0$, define the
probability measure $Q_T(\xi;\cdot)$ by
\begin{eqnarray}\label{def:probab-measure-Q}
Q_{T}(\xi;A)=\frac{1}{T}\int_0^T{\mathcal{P}}_t(\xi;A)\D t,\ \ \ \ \
\ \ \ \ \ A\in{\mathcal{M}},
\end{eqnarray}
where ${\mathcal{P}}_t(\xi;A):={\mathcal{P}}_t\1_A(\xi)$. Recall the $d$-dimensional compensated pure jump process
$\tilde{Y}$ defined as
\eqref{eq:rsddej1}. We have the following result.
\begin{prop}\label{prop:compensatedY}
Suppose that the jump coefficient $g(\cdot)$ satisfies the following linear
growth-type condition:
\begin{eqnarray}\label{cond:g-growth}
|g(t,x,y,\rho)|\leq\ell_g(1+|x|+|y|)h(t,\rho),\ \ \ \ \ \forall\
(t,x,y,\rho)\in\RR_+\times\RR_+^d\times\RR_+^d\times{\mathcal{E}},
\end{eqnarray}
where $h(t,\rho)$ is a nonnegative measurable function satisfying
\begin{eqnarray}\label{cond:h1}
\sup_{t\geq0}\int_{{\mathcal{E}}}h^2(t,\rho)e^{\varrho
h(t,\rho)}\nu(\D\rho)<+\infty,
\end{eqnarray}
for any finite $\varrho>0$. Let $u\geq\tau$. Then for any
$\varepsilon>0$ and $r>0$, there exists a constant
$\delta_{\varepsilon,r}>0$ such that whenever
$\delta\in(0,\delta_{\varepsilon,r}]$,
\begin{eqnarray}\label{eq:modulue-poisson}
{\PP}\left(\sup_{s,t\in[-\tau,0],|s-t|<\delta}\left|\tilde{Y}(u+t)-\tilde{Y}(u+s)\right|\geq
r\right)<\varepsilon.
\end{eqnarray}
\end{prop}

\noindent{\it Proof.}\quad Let $m>0$. Define the stopping time by $
\tau_m=\inf\{t\geq0;\ \|X\|_{[-\tau,t]}>m\}$. Then, $\tau_m\to\infty$ almost surely as $m\to\infty$ due to
Proposition \ref{prop:moment-segment}. Thus we can choose a constant $m>0$ large enough so that
\begin{equation*}\label{eq:stop-beta}
\begin{split}
&{\PP}\left(\sup_{s,t\in[-\tau,0],|s-t|<\delta}\left|\tilde{Y}(u+t)-\tilde{Y}(u+s)\right|\geq r\right)\\
&\quad\leq{\PP}\left(\sup_{s,t\in[-\tau,0],|s-t|<\delta}\left|\tilde{Y}(u+t)-\tilde{Y}(u+s)\right|\geq
r,\tau_{m}>u\right)+{\PP}(\tau_{m}\leq u)\\
&\quad\leq
{\PP}\left(\sup_{s,t\in[-\tau,0],|s-t|<\delta}\left|\tilde{Y}((u+t)\wedge\tau_m)-\tilde{Y}((u+s)\wedge
\tau_m)\right|\geq r\right)+\frac{\varepsilon}{2}.
\end{split}
\end{equation*}
It remains to prove that there exists a constant
$\delta_{\varepsilon,r}>0$ such that whenever
$\delta\in(0,\delta_{\varepsilon,r}]$,
\begin{eqnarray}\label{eq:stop-beta2}
{\PP}\left(\sup_{s,t\in[-\tau,0],|s-t|<\delta}\left|\tilde{Y}((u+t)\wedge\tau_m)-\tilde{Y}((u+s)\wedge
\tau_m)\right|\geq r\right)<\frac{\varepsilon}{2}.
\end{eqnarray}
Using the condition \eqref{cond:g-growth}, the event
\begin{eqnarray*}
\{\tau_m>s\}\subset
\left\{g(s\wedge\tau_m,X((s\wedge\tau_m)),X((s-\tau)\wedge\tau_m),\rho)\leq
\ell_g(1+2m)h(s\wedge\tau_m,\rho)\right\},
\end{eqnarray*}
for any $s>0$. Accordingly, to prove (\ref{eq:stop-beta2}), it
suffices to check that for any $\varepsilon>0$ and $r>0$, there
exists a constant $\delta_{\varepsilon,r}>0$ such that whenever
$\delta\in(0,\delta_{\varepsilon,r}]$,
\begin{eqnarray}\label{eq:stop-beta3}
{\PP}\left(\sup_{s,t\in[-\tau,0],|s-t|<\delta}\left|\tilde{Y}(u+t)-\tilde{Y}(u+s)\right|\geq
r\right)<\frac{\varepsilon}{2},
\end{eqnarray}
under conditions \eqref{cond:g} and \eqref{cond:h1}.

We next take the function $H_\lambda(\cdot)$ on $\RR_+^d$
given by \eqref{eq:fcnH}. Then, by Lemma \ref{lemma:ito-formula},
the process
\begin{eqnarray*}
Y_\lambda^E(t):=\exp\left(H_\lambda(\tilde{Y}(t))-1-\int_0^t\tilde Q^{H_\lambda}(s,X(s),X(s-\tau))\D
s\right),\ \ \ \ t\geq0
\end{eqnarray*}
is a positive ${\mathbb F}$-local martingale, where the function
$\tilde Q^{H_\lambda}(t,x,y)$ with
$(t,x,y)\in\RR_+\times\RR_+^d\times\RR_+^d$ is given by
\begin{eqnarray*}\label{eq:fcnQlambda1}
\tilde Q^{H_\lambda}(t,x,y)=\int_{\mathcal{E}}\left[\exp\left\{{H_\lambda([x+g(t,x,y,\rho)]^+)-H_\lambda(x)}\right\}-1-\left<\nabla
H_\lambda(x),g(t,x,y,\rho)\right>\right]\nu(\D\rho).
\end{eqnarray*}
Further, with conditions \eqref{cond:g} and
\eqref{cond:h1}, we have $\tilde Q^{H_\lambda}(t,x,y)\leq \lambda C$ for all
$(t,x,y)\in\RR_+\times\RR_+^d\times\RR_+^d$, where $C=C(d)>0$ is
some constant depending on the dimension number $d$. Let $u\geq\tau$
and $t\in[-\tau,0]$. For any $r,\lambda>0$, we have
\begin{eqnarray*}
&&{\PP}\left(\frac{\left|\tilde{Y}(u+t)\right|}{\sqrt{|u+t|}}\geq
r\right)={\PP}\left(H_\lambda(\tilde{Y}(u+t))\geq\sqrt{1+\lambda r^2|u+t|}\right)\\
&&\qquad\quad\leq{\PP}\left(1+\int_0^{u+t}\tilde Q^{H_\lambda}(s,X(s),X(s-\tau))\D
s+\log(Y_\lambda^E(u+t))\geq \sqrt{1+\lambda r^2|u+t|}\right)\\
&&\qquad\quad\leq\exp\left[1+\lambda(u+t) C-\sqrt{1+\lambda r^2|u+t|}\right]{\E}\left[\tilde{Y}(u+t)\right]\\
&&\qquad\quad\leq\exp\left[1+ \lambda(u+t) C-\sqrt{\lambda r^2|u+t|}\right],
\end{eqnarray*}
for some constant $C=C(d)>0$. Taking the parameter $\lambda=\beta |u+t|^{-1} r^2$ for some
$\beta>0$, we have
\begin{eqnarray*}
{\PP}\left(\frac{\left|\tilde{Y}(u+t)\right|}{\sqrt{u+t}}\geq
r\right)\leq\exp\left(-r^2\sqrt{\beta}+1+r^2 \beta  C\right).
\end{eqnarray*}
Choose the above constant $\beta>0$ small enough so that
$\beta^*:=\beta(\beta^{-\frac{1}{2}}-C)>0$. Then it holds that
\begin{eqnarray*}
{\PP}\left(\frac{\left|\tilde{Y}(u+t)\right|}{\sqrt{u+t}}\geq
r\right)\leq \exp\left(-\beta^*r^2+1\right).
\end{eqnarray*}
Using the integration by parts, for any $\beta_0\in(0,\beta^*]$, it follows that
\begin{eqnarray*}
{\E}\left[\beta_0\frac{\left|\tilde{Y}(u+t)\right|^2}{|u+t|}\right]<\infty.
\end{eqnarray*}
This further yields that there exists a constant $C>0$ such that
\begin{eqnarray}\label{estimate:module-continuity}
R:=\sup_{s,t\in[-\tau,0],s\neq
t}{\E}\left[\exp\left(C\frac{\left|\tilde{Y}(u+t)-\tilde{Y}(u+s)\right|^2}{|s-t|}\right)\right]<\infty.
\end{eqnarray}
Consider the following positive-valued random variable given by
$$
V=\int_{u-\tau}^u\int_{u-\tau}^u
\exp\left(C\frac{\left|\tilde{Y}(v_1)-\tilde{Y}(v_2)\right|^2}{|v_1-v_2|}\right)\D
v_1\D v_2.
$$
Then ${\E}[V]\leq\tau^2R<\infty$. From Garsia-Rodemich-Rumsey's
Lemma (see \cite[Theorem 2.1.3, p47]{SV}), it follows that
\begin{eqnarray*}
\left|\tilde{Y}(u+t)-\tilde{Y}(u+s)\right|\leq
C\int_0^{|s-t|}\sqrt{\log\left(\frac{V}{v^2}\right)}\D\sqrt{v}.
\end{eqnarray*}
Then there exists a constant $\kappa>0$ small enough such that
\begin{eqnarray}\label{eq:k-Holder}
\left|\tilde{Y}(u+t)-\tilde{Y}(u+s)\right|\leq
C\big(1+\sqrt{\log(V)}\big)|t-s|^\kappa.
\end{eqnarray}
Hence we can arrive at
\begin{equation*}
\begin{split}
&{\PP}\left(\sup_{s,t\in[-\tau,0],|s-t|<\delta}\left|\tilde{Y}(u+t)-\tilde{Y}(u+s)\right|\geq
r\right)
\leq {\PP}\left(C\left[1+\sqrt{\log (V)}\right]\delta^\kappa\geq
r\right)\\
&\qquad \leq{\PP}\left(V\geq
\exp\left[\left|rC^{-1}\delta^{-\kappa}-1\right|^2\right]\right)\leq
{\E}[V]\exp\left[-\left|rC^{-1}\delta^{-\kappa}-1\right|^2\right]\\
&\qquad\leq
\tau^2R\exp\left[-\left|rC^{-1}\delta^{-\kappa}-1\right|^2\right],
\end{split}
\end{equation*}
which yields (\ref{eq:stop-beta3}). \hfill$\Box$

\begin{rem}\label{rem:modulue1}
Let $t\geq0$. Define the ${\mathbb F}$-local
martingale $Z(t)=\int_0^t\sigma(s,X(s),X(s-\tau))\D W(s)$. Kinnally
and Williams \cite{KW} proved that the similar conclusion of
Proposition \ref{prop:compensatedY} holds if the diffusion
coefficient $\sigma(\cdot)$ satisfies the linear growth condition (see Eq.~(23) in \cite{KW}).
\end{rem}

\begin{rem}\label{rem:modulue2}
Let the assumptions in Proposition \ref{prop:compensatedY} hold. Then
for any $\varepsilon>0$ and $r>0$, there exists a constant
$\delta_{\varepsilon,r}>0$ such that whenever
$\delta\in(0,\delta_{\varepsilon,r}]$,
\begin{equation}\label{eq:modulue-poisson2}
\begin{split}
{\PP}\left(\sup_{(t_1<t_2\in[-\tau,0]; |t_1-t_2|<\delta)}
\sup_{t\in[t_1,t_2]}\left[\left|U(u+t)-U(u+t_1)\right|\wedge\left|U(u+t)-U(u+t_2)\right|\right]\geq
r\right)<\varepsilon,
\end{split}
\end{equation}
where $U\in\{Z,\tilde{Y}\}$. In fact, using similar arguments to that
of Proposition \ref{prop:compensatedY}, it suffices to prove the
validity of \eqref{eq:modulue-poisson2} with $U=\tilde{Y}$ under
conditions \eqref{cond:g} and \eqref{cond:h1}. Notice that
\begin{gather}\label{eq:rem-estimate1}
\sup_{(t_1<t_2\in[-\tau,0];
|t_1-t_2|<\delta)}\sup_{t\in[t_1,t_2]}\left[\left|\tilde{Y}(u+t)-\tilde{Y}(u+t_1)\right|
\wedge\left|\tilde{Y}(u+t)-\tilde{Y}(u+t_2)\right|\right]\nonumber\\
\leq2\sup_{(t_1<t_2\in[-\tau,0];
|t_1-t_2|<\delta)}\sup_{t\in[t_1,t_2]}\left|\tilde{Y}(u+t)-\tilde{Y}(u+t_1)\right|.
\end{gather}
Let $t_1<t_2\in[-\tau,0]$ be fixed and satisfies $|t_2-t_2|<\delta$.
Using the estimate \eqref{eq:k-Holder}, for all $t\in[t_1,t_2]$, we have
\begin{eqnarray*}\label{eq:rem-estimate2}
\begin{split}
\left|\tilde{Y}(u+t)-\tilde{Y}(u+t_1)\right|\leq
C\big(1+\sqrt{\log\left(V\right)}\big)|t-t_1|^\kappa \leq
C\big(1+\sqrt{\log\left(V\right)}\big)\delta^\kappa,
\end{split}
\end{eqnarray*}
which implies that
\begin{eqnarray*}\label{eq:rem-estimate3}
\sup_{(t_1<t_2\in[-\tau,0];
|t_1-t_2|<\delta)}\sup_{t\in[t_1,t_2]}\left|\tilde{Y}(u+t)-\tilde{Y}(u+t_1)\right|
\leq C\big(1+\sqrt{\log\left(V\right)}\big)\delta^\kappa.
\end{eqnarray*}
The remaining proof of \eqref{eq:modulue-poisson2} is very similar
to that of Proposition \ref{prop:compensatedY}.
\end{rem}

By virtue of the above auxiliary results on stochastic integrals w.r.t. compensated
Poisson measure $\tilde{N}$ and $n$-dimensional Brownian motion $W$, we have the following
main result.
\begin{theorem}\label{thm:inv-measure}
Let conditions {\sf(A1)} and {\sf (A2)} hold. Suppose that the
jump coefficient $g(\cdot)$ and characteristic measure $\nu(\cdot)$ satisfy
\eqref{cond:g-growth} and \eqref{cond:h1} respectively. The drift
coefficient $b(\cdot)$ satisfy the growth condition:
$|b(t,x,y)|\leq\ell_b(1+|x|^k+|y|^k)$ with $k=1$ or $2$, for all
$(t,x,y)\in\RR_+\times\RR_+^d\times\RR_+^d$, where $\ell_b>0$. Then
the Markov semigroup $({\mathcal{P}}_t;\ t\geq0)$ defined as
\eqref{eq:semigroup} for the segment process of the solution to
Eq.~\eqref{eq:rsddej} admits a unique invariant measure.
\end{theorem}

\noindent{\it Proof.}\quad First we deal with the existence of an
invariant measure for the Markov semigroup $({\mathcal{P}}_t;\ t\geq0)$
defined as \eqref{eq:semigroup}. Let $(t_n;\ n\in{\mathbb{N}})$ be a sequence of times increasing to $+\infty$.
We prove that the sequence of probability measures
$(Q_{t_n}(\xi;\cdot);\ n\in{\mathbb{N}})$ is tight on the Skorohod
space $\{D([-\tau,0];\RR_+^d),{\mathcal{M}}\}$. First the Markov
transition semigroup $({\mathcal{P}}_t(\xi;\cdot);\ t\geq t_0)$ is
Fellerian for some time $t_0\geq0$ using a similar argument to the
non-reflected case documented in Section 3.3 in Reib et al.
\cite{RRG}, any weak limit point is an invariant measure by virtue
of Krylov-Bogulyubov's Theorem (see e.g., \cite{DZ} and \cite{RRG}). We next fix a finite time level $T>0$. Then for any $r>0$,
\begin{eqnarray}\label{eq:QT-estimate0}
Q_T\left(\eta\in D([-\tau,0];\RR_+^d);\ |\eta(0)|>r\right)=\frac{1}{T}\int_0^T\PP\left(\left|X^\xi(t)\right|>r\right)\D t\leq\frac{1}{r^2}\sup_{t\geq0}\E\left[\left|X^\xi(t)\right|^2\right],
\end{eqnarray}
which tends to zero when $r\to\infty$ using Proposition
\ref{prop:moment-segment}. Let $\delta>0$. Define the modulus of
continuity of any function $\eta\in D(I;\RR_+^d)$ with a
subinterval $I$ of $\RR$:
\begin{eqnarray*}
w(\eta;\delta)&:=&\sup_{(s,t\in I;\ |s-t|<\delta)}|\eta(s)-\eta(t)|,\ \ \ \ \ \ {\rm and}\nonumber\\
w^*(\eta;\delta)&:=&\sup_{(t_1<t_2\in I;\ t_2-t_1<\delta)}\sup_{t_1\leq t\leq t_2}\left[|\eta(t)-\eta(t_1)|\wedge|\eta(t)-\eta(t_2)|\right].
\end{eqnarray*}
Using the solution representation of Skorohod problem (see e.g., Asmussen~\cite{AS03}), the regulator $K=(K^i(t);\ t\geq0)_{i=1,\ldots,d}$ admits
the representation given by
\begin{eqnarray}\label{eq:skorohod-solution}
K(t)=\sup_{0\leq s\leq t}\left[\Gamma(X)(s)\right]^-,\ \ \ \ \ \ t\geq0,
\end{eqnarray}
where the $d$-dimensional process $\Gamma(X)=(\Gamma(X)(t);\ t\geq0)$ is defined as
\begin{eqnarray}\label{eq:Gamma}
\Gamma(X)(t)&=&\xi(0) + \int_0^t b(s,X(s),X(s-\tau))\D s + \int_0^t\sigma(s,X(s),X(s-\tau))\D W(s)\nonumber\\
&&+\int_0^t\int_{{\mathcal{E}}}g(s,X(s-),X((s-\tau)-),\rho)\tilde{N}(\D\rho,\D
s).
\end{eqnarray}
Then for any $u\geq\tau$, we have
\begin{eqnarray*}
&&\sup_{(s<t\in[-\tau,0];\ |s-t|<\delta)}\left|X(u+t)-X(u+s)\right|\nonumber\\
&&\qquad=\sup_{(s<t\in[-\tau,0];\ |s-t|<\delta)}\left|\Gamma(X)(u+t)-\Gamma(X)(u+s)+K(u+t)-K(u+s)\right|\nonumber\\
&&\qquad\leq2\sup_{(s<t\in[-\tau,0];\ |s-t|<\delta)}\left|\Gamma(X)(u+t)-\Gamma(X)(u+s)\right|\nonumber\\
&&\qquad\leq2\sup_{(s<t\in[-\tau,0];\ |s-t|<\delta)}\int_{u+s}^{u+t}\left|b(v,X(v),X(v-\tau))\right|\D v+2\sup_{(s<t\in[-\tau,0];\ |s-t|<\delta)}\left|Z(u+t)-Z(u+s)\right|\nonumber\\
&&\qquad\qquad+2\sup_{(s<t\in[-\tau,0];\ |s-t|<\delta)}\left|\tilde{Y}(u+t)-\tilde{Y}(u+s)\right|.
\end{eqnarray*}
Here the stochastic integral processes $Z$ and
$\tilde{Y}$ are defined as in Remark \ref{rem:modulue1} and Remark \ref{rem:modulue2} respectively.
Hence for any $r>0$ and $k=1$ or $2$,
\begin{eqnarray}\label{eq:prob-modu}
\PP\left(w(X_u;\delta)\geq r\right)&\leq&\PP\left(\sup_{(s<t\in[-\tau,0];\ |s-t|<\delta)}
\int_{u+s}^{u+t}[1+|X(v)|^k+|X(v-\tau)|^k]\D v\geq\frac{r}{6\ell_b}\right)\nonumber\\
&&+\PP\left(\sup_{(s<t\in[-\tau,0];\ |s-t|<\delta)}\left|Z(u+t)-Z(u+s)\right|\geq\frac{r}{6}\right)\nonumber\\
&&+\PP\left(\sup_{(s<t\in[-\tau,0];\ |s-t|<\delta)}\left|\tilde{Y}(u+t)-\tilde{Y}(u+s)\right|\geq\frac{r}{6}\right)\nonumber\\
&=:&F_1(u,\delta;r)+F_2(u,\delta;r)+F_3(u,\delta;r).
\end{eqnarray}
Take the following inequalities into account
\begin{gather*}
\sup_{u\geq\tau}\PP\left(\|X\|_{[u-2\tau,u]}^k>r\right)\leq\left\{\begin{array}{cl}
\frac{1}{r^2}\sup_{t\geq0}\E\left[\|X_t\|_{[-\tau,0]}^2\right], & k=1,\\ \\
\frac{1}{r}\sup_{t\geq0}\E\left[\|X_t\|_{[-\tau,0]}^2\right], & k=2,
\end{array}\right.\nonumber\\
F_1(u,\delta;r)\leq\PP\left(\delta[C_1+C_2\|X\|_{[u-2\tau,u]}^k]\geq\frac{r}{6\ell_b}\right),\ \ \ {\rm for}\ k=1\ {\rm or}\ 2.
\end{gather*}
By virtue of \eqref{eq:moment-segment} in Proposition
\ref{prop:moment-segment}, for any $\varepsilon,r>0$, there exists a
constant $\delta_{\varepsilon,r}^1>0$ such that
$\sup_{u\geq\tau}F_1(u,\delta;r)<\frac{\varepsilon}{6}$
whenever $\delta\in(0,\delta_{\varepsilon,r}^1]$. It follows from
Proposition \ref{prop:compensatedY} and Remark \ref{rem:modulue1}
that there exist constants $\delta_{\varepsilon,r}^2>0$ and
$\delta_{\varepsilon,r}^3>0$ such that
$\sup_{u\geq\tau}F_2(u,\delta;r)<\frac{\varepsilon}{6}$ when
$\delta\in(0,\delta_{\varepsilon,r}^2]$ and
$\sup_{u\geq\tau}F_3(u,\delta;r)<\frac{\varepsilon}{6}$ for
$\delta\in(0,\delta_{\varepsilon,r}^3]$ respectively. Finally, we
obtain
\begin{eqnarray}\label{eq:w-estimate}
\sup_{u\geq\tau}\PP\left(w(X_u;\delta)\geq r\right)<\frac{\varepsilon}{2},
\end{eqnarray}
whenever $\delta\in(0,\wedge_{i=1}^3\delta_{\varepsilon,r}^i]$. For the above $\delta$, we further conclude that for all
$T>\frac{2\tau}{\varepsilon}\vee\tau$,
\begin{eqnarray}\label{eq:Q-estimate1}
Q_T\left(\eta\in D([-\tau,0];\RR_+^d);\ w(\eta;\delta)\geq r\right)\leq\frac{\tau}{T}
+\frac{1}{T}\int_\tau^T\PP\left(w(X_u;\delta)\geq r\right)\D u<\varepsilon.
\end{eqnarray}
For any $\varepsilon,r>0$, we can conclude that there exists a
constant $\delta^*_{\varepsilon,r}>0$ so that for all
$\delta\in(0,\delta^*_{\varepsilon,r}]$, $Q_T\left(\eta\in
D([-\tau,0];\RR_+^d);\ w^*(\eta;\delta)\geq r\right)<\varepsilon$
by employing Remark \ref{rem:modulue2}. Together with
\eqref{eq:QT-estimate0} and \eqref{eq:Q-estimate1}, the sequence of
probability measures $(Q_{t_n}(\xi;\cdot);\ n\in{\mathbb{N}})$ is
tight on the Skorohod space $\{D([-\tau,0];\RR_+^d),{\mathcal{M}}\}$
using Theorem 6.6 in Liptser and Shiryayev \cite{LS}.

Next we check the uniqueness of invariant measures under the assumption {\sf(A2)}. Let $X^\xi=(X^\xi(t);\ t\geq-\tau)$ and
$X^\eta=(X^\eta(t);\ t\geq-\tau)$ be two strong solutions to Eq.~\eqref{eq:rsddej} with respect initial datum $\xi,\eta\in D([-\tau,0];\RR_+^d)$. Then for
$t\geq0$,
\begin{eqnarray*}
&& X^\xi(t)-X^\eta(t)=\xi(0)-\eta(0) + \int_0^t [b(s,X^\xi(s),X^\xi(s-\tau))-b(s,X^\eta(s),X^\eta(s-\tau))]\D s\nonumber\\
&&\qquad+\int_0^t [\sigma(s,X^\xi(s),X^\xi(s-\tau))-\sigma(s,X^\eta(s),X^\eta(s-\tau))]\D W(s)\nonumber\\
&&\qquad+\int_0^t\int_{{\mathcal{E}}}[g(s,X^\xi(s-),X^\xi((s-\tau)-),\rho)-g(s,X^\eta(s-),X^\eta((s-\tau)-),\rho)]\tilde{N}(\D\rho,\D s)\nonumber\\
&&\qquad+K^\xi(t)-K^\eta(t),
\end{eqnarray*}
where $K^\xi=(K^\xi(t);\ t\geq0)$ and $K^\eta=(K^\eta(t);\ t\geq0)$ denotes the regulators for the solutions $X^\xi$ and $X^\eta$
respectively. For convenience, we let
$b^j(t):=b(t,X^j(t),X^j(t-\tau))$,
$\sigma^j(t):=\sigma(t,X^j(t),X^j(t-\tau))$ and
$g^j(t,\rho):=g(t,X^j(t),X^j(t-\tau),\rho)$ with $j\in\{\xi,\eta\}$.
Using the It\^o's formula \eqref{eq:ito-formula}, we arrive at
\begin{eqnarray}\label{eq:ito-diff}
\left|X^\xi(t)-X^\eta(t)\right|^2&=&\left|\xi(0)-\eta(0)\right|^2 + 2\int_0^t \left<(X^\xi-X^\eta)(s),(b^\xi-b^\eta)(s)\right>\D s\nonumber\\
&&+2\int_0^t \left<(X^\xi-X^\eta)(s),(\sigma^\xi-\sigma^\eta)(s)\D W(s)\right>+\int_0^t\left\|(\sigma^\xi-\sigma^\eta)(s)\right\|^2\D s\nonumber\\
&&+2\int_0^t\left<(X^\xi-X^\eta)(s),\D K^{\xi,c}(s)-\D K^{\eta,c}(s)\right>\nonumber\\
&&+\sum_{0<s\leq t}\Big[|(X^\xi-X^\eta)(s-)+\Delta(X^\xi-X^\eta)(s)|^2-|(X^\xi-X^\eta)(s-)|^2\nonumber\\
&&\qquad\qquad-2\left<(X^\xi-X^\eta)(s-),\Delta(Y^\xi-Y^\eta)(s)\right>\bigg].
\end{eqnarray}
Here the pure jump processes
$Y^j(t)=\int_0^t\int_{{\mathcal{E}}}g^j(s-,\rho)N(\D\rho,\D s)$ and
$K^{j,c}(t)$ corresponds to the continuous counterpart of $K^j(t)$
with $j\in\{\xi,\eta\}$. It follows from \eqref{eq:jump-Xi} that
\begin{eqnarray*}
\Delta(X^\xi-X^\eta)(t)=\varphi(X^\xi(t-),\Delta
Y^\xi(t))-\varphi(X^\xi(t-),\Delta Y^\xi(t)),
\end{eqnarray*}
where the function $\varphi(x,y)$ is defined in \eqref{eq:jump-Xi}
with $x,y\in\RR_+^d$. Then it holds that
\begin{eqnarray*}
(X^\xi-X^\eta)(t-)+\Delta(X^\xi-X^\eta)(t)&=&[X^\xi(t-)+\varphi(X^\xi(t-),\Delta Y^\xi(t))]-[X^\eta(t-)+\varphi(X^\eta(t-),\Delta Y^\eta(t))]\nonumber\\
&=&[X^\xi(t-)+\Delta Y^\xi(t)]^+-[X^\eta(t-)+\Delta Y^\eta(t)]^+.
\end{eqnarray*}
Thus Eq.~\eqref{eq:ito-diff} becomes
\begin{eqnarray*}
&&\left|X^\xi(t)-X^\eta(t)\right|^2=\left|\xi(0)-\eta(0)\right|^2 + 2\int_0^t \left<(X^\xi-X^\eta)(s),(b^\xi-b^\eta)(s)\right>\D s\nonumber\\
&&\qquad+2\int_0^t \left<(X^\xi-X^\eta)(s),(\sigma^\xi-\sigma^\eta)(s)\D W(s)\right>+\int_0^t\|(\sigma^\xi-\sigma^\eta)(s)\|^2\D s\nonumber\\
&&\qquad+2\int_0^t\left<(X^\xi-X^\eta)(s),\D K^{\xi,c}(s)-\D K^{\eta,c}(s)\right>\nonumber\\
&&\qquad+\int_0^t\int_{{\mathcal{E}}}\Big[\left|[X^\xi(s-)+g^\xi(s-,\rho)]^+-[X^\eta(s-)+ g^\eta(s-,\rho)]^+\right|^2-|(X^\xi-X^\eta)(s-)|^2\nonumber\\
&&\qquad\qquad-2\left<(X^\xi-X^\eta)(s-),(g^\xi-g^{\eta})(s-,\rho)\right>\Big]N(\D\rho,\D
s).
\end{eqnarray*}
Since the function $x\to [x]^+$ is Lipschitzian continuous, it holds that
\[
\left|[X^\xi(t-)+g^\xi(t,\rho)]^+-[X^\eta(t-)+
g^\eta(t,\rho)]^+\right|^2\leq
\left|(X^\xi-X^\eta)(t-)+(g^\xi-g^{\eta})(t,\rho)\right|^2.
\]
Let $\lambda>0$. Take the condition {\sf(A2)} into account, we have
\begin{eqnarray*}
\E\left[e^{\lambda t}\left|X^\xi(t)-X^\eta(t)\right|^2\right]&\leq&\E\left[\left|\xi(0)-\eta(0)\right|^2\right] + 2\E\left[\int_0^te^{\lambda s} \left<(X^\xi-X^\eta)(s),(b^\xi-b^\eta)(s)\right>\D s\right]\nonumber\\
&&+\E\left[\int_0^te^{\lambda s}\|(\sigma^\xi-\sigma^\eta)(s)\|^2\D s\right]+\lambda\E\left[\int_0^te^{\lambda s}\left|X^\xi(s)-X^\eta(s)\right|^2\D s\right]\nonumber\\
&&+\E\left[\int_0^t\int_{{\mathcal{E}}}e^{\lambda s}|(g^\xi-g^{\eta})(s,\rho)|^2\nu(\D\rho)\D s\right]\nonumber\\
&\leq&\E\left[|\xi(0)-\eta(0)|^2\right] +e^{\lambda\tau}\E\left[\int_{-\tau}^0e^{\lambda v}|\xi(v)-\eta(v)|^2\D v\right]\nonumber\\
&&+(\lambda-\alpha_1+\alpha_2e^{\lambda\tau})\E\left[\int_0^te^{\lambda s}|(X^\xi-X^\eta)(s)|^2\D s\right],
\end{eqnarray*}
where we used the following estimate of regulators $K^\xi$ and $K^\eta$, for $t\geq0$,
\begin{eqnarray*}
&&\int_0^te^{\lambda s}\left<(X^\xi-X^\eta)(s),\D K^{\xi,c}(s)-\D K^{\eta,c}(s)\right>\nonumber\\
&&\qquad=\int_0^te^{\lambda s}\left<X^\xi(s),\D
K^{\xi,c}(s)\right>-\int_0^te^{\lambda s}\left<X^\eta(s),\D
K^{\xi,c}(s)\right>
-\int_0^te^{\lambda s}\left<X^\xi(s),\D K^{\eta,c}(s)\right>\nonumber\\
&&\qquad\qquad+\int_0^te^{\lambda s}\left<X^\eta(s),\D K^{\eta,c}(s)\right>\nonumber\\
&&\qquad=-\int_0^te^{\lambda s}\left<X^\eta(s),\D K^{\xi,c}(s)\right>-\int_0^te^{\lambda s}\left<X^\xi(s),\D K^{\eta,c}(s)\right>\leq0.
\end{eqnarray*}
Here we have used the fact that, for $j\in\{\xi,\eta\}$, it holds using the support property \eqref{eq:Kcsupport} that
\begin{eqnarray*}
0\leq \int_0^t e^{\lambda s}\left<X^j(s),\D K^{j,c}(s)\right>\leq e^{\lambda t}\int_0^t \left<X^j(s),\D K^{j,c}(s)\right>=0,\ \ \forall\ t\geq0.
\end{eqnarray*}
Recall the positive constant $\lambda^*$ satisfying $\lambda^*-\alpha_1+\alpha_2e^{\lambda^*\tau}=0$. Then for all $t\geq0$,
\begin{eqnarray*}
\E\left[\left|(X^\xi-X^\eta)(t)\right|^2\right]\leq e^{-\lambda^*
t}\left\{\E\left[|\xi(0)-\eta(0)|^2\right]
+e^{\lambda^*\tau}\E\left[\int_{-\tau}^0e^{\lambda^*
v}|\xi(v)-\eta(v)|^2\D v\right]\right\},
\end{eqnarray*}
which shows that
\begin{eqnarray}\label{eq:diff-0}
\lim_{t\to+\infty}\E\left[\left|(X^\xi-X^\eta)(t)\right|^2\right]=0.
\end{eqnarray}
Using the similar proof to that of \eqref{eq:moment-single}, we have
\begin{eqnarray*}
\E\left[\left\|X_t^\xi-X_t^\eta\right\|_{[-\tau,0]}^2\right]\leq
4\E\left[\left|(X^\xi-X^\eta)(t-\tau)\right|^2\right]+
C\int_{t-\tau}^t\E\left[\left|(X^\xi-X^\eta)(s-\tau)\right|^2\right]\D
s,
\end{eqnarray*}
for some constant $C>0$ which is independent of time $t$. The above estimate and Gronwall's lemma lead to
$$\lim_{t\to+\infty}\E\left[\left\|X_t^\xi-X_t^\eta\right\|_{[-\tau,0]}^2\right]=0.$$
Thus we complete the proof of the uniqueness of invariant measures.
\hfill$\Box$

\section{Local Time of Solutions}
\label{sec:local-time}

As we have pointed out in Section~\ref{sec:intro}, the sizes of the
jumps for the regulator $K=(K^i(t);\ t\geq0)_{i=1,\ldots,d}$ can be
identified by Eq.~\eqref{eq:jumpK}. In this section, we will establish a
relationship between the continuous counterpart of the regulator $K$
and the local time of the solution process.

For $i=1,\dots,d$, let $L^i=(L^i(t);\ t\geq0)$ be the local time
process for the $i$th-element of the strong solution process $X$ at
point $0$. Moreover, if $\sum_{0<s\leq t}|\Delta Y^i(s)|<+\infty$
for all $t>0$ (the pure jump process $Y^i$ is defined as
\eqref{eq:Y-pure-jump}), then the local time has the following limit
representation (see Protter \cite{Protter}):
\begin{eqnarray}\label{eq:local-times}
L^i(t)&=&\lim_{\varepsilon\downarrow0}\frac{1}{\varepsilon}\int_0^t\1_{[0,\varepsilon)}(X^i(s))\D\left<X^{i,c},X^{i,c}\right>_s\nonumber\\
&=&\lim_{\varepsilon\downarrow0}\frac{1}{\varepsilon}\sum_{j=1}^n\int_0^t\1_{[0,\varepsilon)}(X^i(s))\sigma_{ij}^2(s,X(s),X(s-\tau))\D s,
\end{eqnarray}
where $n\in{\mathbb{N}}$ is the dimension of the Brownain motion $W$. Further define the function by, for $i=1,\dots,d$,
\begin{eqnarray}\label{eq:fcnhatb}
\hat{b}_i(t,x,y):=b_i(t,x,y)-\int_{\mathcal{E}}g_i(t,x,y,\rho)\nu(\D\rho),\
\ \ \ \ (t,x,y)\in\RR_+\times\RR_+^d\times\RR_+^d.
\end{eqnarray}
Then we have the following relationship between the local time $L^i$ and the $i$th-regulator $K^i$.
\begin{prop}\label{prop:local-time}
For $i=1,\dots,d$, it holds that
\begin{eqnarray}\label{eq:local-timeK}
\frac{1}{2}L^i(t)=\int_0^t\1_{\{X^i(s)=0\}}\hat{b}_i(s,X(s),X(s-\tau))\D
s+K^{i,c}(t),\ \ \ \ \ \ \ t\geq0,
\end{eqnarray}
where $K^{i,c}(t)$ is the continuous counterpart of the
$i$th-regulator $K^{i}(t)$. Moreover, if there exists a positive Borel
measurable function $f_i(t,x)$ on $(t,x)\in\RR_+\times\RR_+^d$ and a
positive Borel measurable function $l_i(x_i)$ on $x\in\RR_+$ such
that
\begin{eqnarray}\label{eq:cond-localtime}
\sum_{j=1}^n\sigma_{ij}^2(t,x,y)\geq l_i(x_i),\ \ \ \ \
\big|\hat{b}_i(t,x,y)\big|\leq f_i(t,x),
\end{eqnarray}
with $(t,x,y)\in\RR_+\times\RR_+\times\RR_+$, then it holds that
\begin{eqnarray}\label{eq:local-timeK2}
\frac{1}{2}L^i(t)=K^{i,c}(t),\ \ \ \ \ \ \ t\geq0.
\end{eqnarray}
\end{prop}

\noindent{\it Proof.}\quad In terms of Eq.~\eqref{eq:rsddej}, the $i$th-element of
the solution process $X$ is given by
\begin{eqnarray*}
X^{i}(t)&=&\xi^{i}(0) + \int_0^tb_i(s,X(s),X(s-\tau))\D
s+\sum_{j=1}^n\int_0^t\sigma_{ij}(s,X(s),X(s-\tau))\D
W^j(s)\nonumber\\
&&+\int_0^t\int_{\mathcal{E}}g_i(s,X(s-),X((s-\tau)-),\rho)\tilde{N}(\D\rho,\D
s) + K^i(t)\nonumber\\
&\geq&0,\ \ \ \ \ \ \ \ \ \ {\rm on}\ t\geq0,\nonumber\\
X^i(t)&=&\xi^{i}(t),\ \ \ \ \ \ {\rm on}\ t\in[-\tau,0].
\end{eqnarray*}
Since the process $(X^i(t))_{t\geq0}$ is a r.c.l.l. semimartingale, using Tanaka's formula (see Protter
\cite{Protter}), for $t\geq0$, we have
\begin{eqnarray*}
X^i(t)&=&\xi^i(0) + \int_0^t\1_{\{X^i(s-)>0\}}\D X^i(s) +
\sum_{0<s\leq t}\1_{\{X^i(s-)=0\}}X^i(s)
+\frac{1}{2}L^i(t)\nonumber\\
&=&X^i(t)-\int_0^t\1_{\{X^i(s-)=0\}}\D X^i(s)+ \sum_{0<s\leq
t}\1_{\{X^i(s-)=0\}}X^i(s) +\frac{1}{2}L^i(t).
\end{eqnarray*}
As a consequence
\begin{eqnarray*}
\frac{1}{2}L^i(t)&=&\int_0^t\1_{\{X^i(s-)=0\}}\D
X^i(s)-\sum_{0<s\leq t}\1_{\{X^i(s-)=0\}}X^i(s)\nonumber\\
&=&\int_0^t\1_{\{X^i(s-)=0\}}\D
X^i(s)-\sum_{0<s\leq t}\1_{\{X^i(s-)=0\}}[X^i(s)-X^i(s-)]\nonumber\\
&=&\int_0^t\1_{\{X^i(s-)=0\}}\D
X^i(s)-\sum_{0<s\leq t}\1_{\{X^i(s-)=0\}}\Delta X^i(s)\nonumber\\
&=&\int_0^t\1_{\{X^i(s-)=0\}}\D X^{i,c}(s),
\end{eqnarray*}
where $X^{i,c}(t)$ corresponds to the continuous part of $X^i(t)$. For $i=1,\dots,d$,
define the process  by
\[
M_i^W(t)=\sum_{j=1}^n\int_0^t\1_{\{X^i(s-)=0\}}\sigma_{ij}(s,X(s),X(s-\tau))\D
W^j(s),\ \ \ \ \ \ t\geq0.
\]
Then for $t\geq0$, we have
\begin{eqnarray}\label{eq:LK1}
\frac{1}{2}L^i(t)=\int_0^t\1_{\{X^i(s-)=0\}}\hat{b}_i(s,X(s),X(s-\tau))\D
s+M^W_i(t)+K^{i,c}(t),
\end{eqnarray}
which implies that the process $M_i^W=(M_i^W(t);\ t\geq0)$ is of
finite variation, since $K^{i,c}(t)$ is non-decreasing w.r.t. time $t\geq0$. Note that $M_i^W$ is also an $\mathbb{F}$-local martingale. Then it must hold that $M_i^W(t)=M_i^W(0)=0$ for $t\geq0$. Using
\eqref{eq:LK1} again, it follows that
\begin{eqnarray*}\label{eq:LK2}
\frac{1}{2}L^i(t)=\int_0^t\1_{\{X^i(s-)=0\}}\hat{b}_i(s,X(s),X(s-\tau))\D
s+K^{i,c}(t),
\end{eqnarray*}
where the function $\hat{b}_i$ is defined as
\eqref{eq:fcnhatb} with $i=1,\dots,d$. This shows the validity of
\eqref{eq:local-timeK}.

Next we verify the validity of \eqref{eq:local-timeK2} under the
condition \eqref{eq:cond-localtime}. As a matter of fact, as a
simple consequence of the occupation time formula (see Exercise
(1.15) in Revuz and Yor \cite{RM91}), we have for $t\geq0$,
\begin{align*}
&\int_0^t\1_{\{X^i(s)=0\}}\big|\hat{b}_i(s,X(s),X(s-\tau))\big|\D\left<X^{i,c},X^{i,c}\right>_s
\leq\int_0^t\1_{\{X^i(s)=0\}}f_i(s,X(s))\D\left<X^{i,c},X^{i,c}\right>_s\nonumber\\
&\qquad=\int_{0}^\infty\left(\int_0^t\1_{\{a=0\}}f_i(s,a)\D
L^{i,a}(s)\right)\D a=\int_{0}^\infty\left(\int_0^t\1_{\{a=0\}}f_i(s,0)\D
L^{i}(s)\right)\D a\nonumber\\
&\qquad=\left(\int_0^tf_i(s,0)\D L^{i}(s)\right)\int_{0}^\infty\1_{\{a=0\}}\D a=0,
\end{align*}
where the nonnegative process $L^{i,a}=(L^{i,a}(t);\ t\geq0)$ denotes the local time of the ith-element $X^i$ of the solution
process $X$ at point $a\geq0$. Thus we obtain
\begin{eqnarray*}
\int_0^t\1_{\{X^i(s)=0\}}\big|\hat{b}_i(s,X(s),X(s-\tau))\big|l_i(X^i(s))\D
s=0,\ \ \  \ \forall\ t\geq0.
\end{eqnarray*}
This yields that $\int_0^t\1_{\{X^i(s)=0\}}\hat{b}_i(s,X(s),X(s-\tau))\D s=0$ for all $t\geq0$ which proves the validity of the equality \eqref{eq:local-timeK2}.

\begin{rem}\label{rem:example-localtime}
We present an illustrative example for the condition \eqref{eq:cond-localtime} in the case of the dimension number
$d=n=1$. Let $(t,x,y,\rho)\in\RR_+^3\times{\mathcal{E}}$. We take the drift coefficient $b(t,x,y)=-\gamma(t)x+\theta_1(t)y$, the diffusion coefficient
$\sigma(t,x,y)=l_1(x)+l_2(t,y)$ and the jump coefficient $g(t,x,y,\rho)=(\ell_g(t)x+\theta_2(t)y)h(t,\rho)$,
where $\gamma(t),l_1(x),l_2(t,y),h(t,\rho)>0$ and $\theta_1(t),\theta_2(t),\ell_g(t)\in\RR$. For all $t\geq0$, assume that $\ell_h(t):=\int_{{\mathcal{E}}}
h(t,\rho)\nu(\D\rho)$ is finite and the positive functions $x\to l_1(x)$ and $y\to l(t,y)$ are Lip-continuous with respect Lip-constants $\ell_1,\ell_2>0$. We take $\theta_1(t)=\theta_2(t)\ell_h(t)$ and choose appropriate set of parameters $(\gamma(t),\theta_2(t),\ell_g(t),\ell_h(t),\ell_1,\ell_2)$ such that $(b,\sigma,g)$ satisfies conditions {\sf(A1)} and {\sf(A2)} (see the illustrative example presented in Section \ref{sec:intro}). In this case, we also have $|\hat{b}(t,x,y)|=|b(t,x,y)-\int_{{\mathcal{E}}}g(t,x,y,\rho)\nu(\D\rho)|
\leq|\gamma(t)+\ell_h(t)+\ell_g(t)|x$ and $\sigma^2(t,x,y)\geq\ell_1^2(x)$ with $(t,x,y)\in\RR_+^3$. Thus the condition \eqref{eq:cond-localtime}
holds.
\end{rem}
\begin{rem}\label{rem:local-time-asymoptic}
\begin{itemize}
  \item[{\sf(1)}] Using the equalities \eqref{eq:local-timeK} and
\eqref{eq:jumpK}, we can characterize the $i$th-regulator $K^i$ in Eq.~\eqref{eq:rsddej} by the following way:
\begin{align}\label{eq:Ki}
K^i(t)&=\frac{1}{2}L^i(t)-\int_0^t\1_{\{X^i(s)=0\}}\hat{b}_i(s,X(s),X(s-\tau))\D
s\nonumber\\
&\quad+\sum_{0<s\leq
t}\left[\int_{\mathcal{E}}g_i(s,X(s-),X((s-\tau)-),\rho)N(\D\rho,\{s\}) + X^i(s-)\right]^-,
\end{align}
for all $t>0$. If the condition \eqref{eq:cond-localtime} holds,
then
\begin{eqnarray}\label{eq:Ki2}
K^i(t)=\frac{1}{2}L^i(t)+\sum_{0<s\leq
t}\left[\int_{\mathcal{E}}g_i(s,X(s-),X((s-\tau)-),\rho)N(\D\rho,\{s\}) + X^i(s-)\right]^-.
\end{eqnarray}
\item[{\sf(2)}] If the jump coefficient $g(\cdot)$ is nonnegative, then the
regulator $K=(K^i(t);\ t\geq0)_{i=1,\ldots,d}$ has a continuous path
modification by \eqref{eq:jumpK}. In this case, assume the condition
\eqref{eq:cond-localtime} is satisfied (see Remark
\ref{rem:example-localtime}), then by Proposition
\ref{prop:local-time}, it holds that
\begin{align*}
L^i(t)&=2X^i(t)-2\xi^i(0)-2\int_0^tb_i(s,X(s),X(s-\tau))\D s-2\sum_{j=1}^n\int_0^t\sigma_{ij}(s,X(s),X(s-\tau),\rho)\D W^j(s)\nonumber\\
&\quad+2\int_0^t\int_{{\mathcal{E}}}g_i(s,X(s-),X(s-\tau)-)\tilde{N}(\D\rho,\D s).
\end{align*}
Hence in the stationary setting, we have
\begin{eqnarray*}
\frac{\E\left[L^i(t)\right]}{t}=-\frac{2}{t}\int_0^t\E\left[b_i(s,\xi(0),\xi(-\tau))\right]\D s,\ \ \ \ \ \ t>0.
\end{eqnarray*}
If the drift function $b(\cdot)$ is independent of time $t$, then
\begin{eqnarray*}
\lim_{t\to+\infty}\frac{\E\left[L^i(t)\right]}{t}=-2\E\left[b_i(\xi(0),\xi(-\tau))\right].
\end{eqnarray*}
The above quantity is usually called the loss rate in the reflected dynamics (see e.g., Asmussen \cite{AS03}).
\end{itemize}
\end{rem}

\section*{Acknowledgments}
The authors would like to thank the associate editor and the reviewer for the valuable comments and suggestions to improve the paper greatly. The research of the first author is partially supported by NCET-12-0914 and NSF of China (grant no. 11471254).

\end{document}